\documentclass[10pt,a4paper]{article}
\usepackage[nice]{nicefrac}
\usepackage{cite}
\usepackage{amscd}
\usepackage[latin1]{inputenc}
\usepackage{amsmath}
\usepackage{amsfonts}
\usepackage{amssymb}
\usepackage{color}
\usepackage{float}
\usepackage{amsthm}
\usepackage{graphicx}
\usepackage{listings}
\usepackage{hyperref}
\usepackage{mathtools}

\newtheorem*{theorem*}{Theorem}
\newtheorem{theorem}{Theorem}

\def\Qed{\hfill\raisebox{.6ex}{\framebox[2.5mm]{}}\\[.15in]}

\def\m{\mathbb}

\makeatletter
\newcommand{\xdashrightarrow}[2][]{\ext@arrow 0359\rightarrowfill@@{#1}{#2}}
\newcommand{\xdashleftarrow}[2][]{\ext@arrow 3095\leftarrowfill@@{#1}{#2}}
\newcommand{\xdashleftrightarrow}[2][]{\ext@arrow 3359\leftrightarrowfill@@{#1}{#2}}
\def\rightarrowfill@@{\arrowfill@@\relax\relbar\rightarrow}
\def\leftarrowfill@@{\arrowfill@@\leftarrow\relbar\relax}
\def\leftrightarrowfill@@{\arrowfill@@\leftarrow\relbar\rightarrow}
\def\arrowfill@@#1#2#3#4{%
  $\m@th\thickmuskip0mu\medmuskip\thickmuskip\thinmuskip\thickmuskip
   \relax#4#1
   \xleaders\hbox{$#4#2$}\hfill
   #3$%
}
\makeatother

\begin{document}
\date{}
\title{Surfaces with canonical map of maximum degree}
\author{Carlos Rito}
\maketitle

\begin{abstract}
We use the Borisov-Keum equations of a fake projective plane and the Borisov-Yeung equations of the
Cartwright-Steger surface to show the existence of a regular surface with canonical map of degree 36
and of an irregular surface with canonical map of degree 27.
As a by-product, we get equations (over a finite field)
for the $\mathbb Z/3$-invariant fibres of the Albanese fibration of the Cartwright-Steger surface and show that they are smooth.

\noindent 2010 MSC: 14J29, 14Q05, 14Q10.

\end{abstract}

{\bf Keywords:}  Surface of general type, Canonical map, Ball-quotient surface.

\section{Introduction}

Let $S$ be a smooth minimal surface of general type with geometric genus $p_g\geq 3,$ irregularity $q$ and
self-intersection of the canonical divisor $K^2.$
Denote by $\phi=\phi_S$ the canonical map of $S$ and let $d:=\deg(\phi).$
Beauville \cite{Be} has proved that, if $d$ is finite, then
$$d\leq 36\ \ {\rm if}\ \ q=0,\ \ \ \ \ d\leq 27\ \ {\rm if}\ \ q>0.$$

Only recently examples with $d>16$ have been given, see \cite{GPR}, \cite{Ri} for $d=24$ and \cite{GPR} for $d=32.$ 
It follows from Beauville's proof that the limit cases $d=36,q=0$ and $d=27,q>0$ can only occur for surfaces with invariants
\begin{equation}\label{inv}
p_g=3,q=0,K^2=36\ \ \ \ {\rm and}\ \ \ \ p_g=3,q=1,K^2=27,
\end{equation}
 respectively. These satisfy $K^2=9\chi,$ hence are ball-quotient surfaces.
 
Surfaces of general type with invariants $K^2=9\chi=9$ and $p_g=0$ (thus $q=0$) are the so-called {\em fake projective planes}.
There are $50$ pairs of complex-conjugated such surfaces,
according to the results of Prasad and Yeung \cite{PY1}, \cite{PY2}, and Cartwright and Steger \cite{CS},
who have also found the unique known example of a surface with invariants
$K^2=9,p_g=q=1$ (the so-called Cartwright-Steger surface). 
 
The only surfaces available in the literature with invariants (\ref{inv}) are certain \'etale coverings of fake projective planes and
of the Cartwright-Steger surface.
In order to prove that their canonical map is of maximum degree, it suffices to show that the canonical system is free from base points.
Since these surfaces are given by uniformization only, this is a hard task. But recently two papers appeared, Borisov-Keum \cite{BK} and
Borisov-Yeung \cite{BY}, giving equations for a (pair of) fake projective plane $Z$ and for the Catwright-Steger surface $S,$ both embedded
in $\mathbb P^9$ by the bicanonical map.

For a long time people have searched for a more explicit construction of such surfaces, so these results were received with enthusiasm.
But the equations are not nice, in the sense that computations are hard even for powerful computers. In this paper we show that we can actually
prove results using their equations, namely we prove that:
\begin{theorem}
Let $Z$ be the above fake projective plane and $S$ be the Cartwright-Steger surface.
Denote by $\phi_X$ the canonical map of $X$. We have that:
\newline
There is an \'etale $(\mathbb Z/2)^2$-covering $\tilde Z\rightarrow Z$ such that $\deg(\phi_{\tilde Z})=36;$
\newline
There is an \'etale $\mathbb Z/3$-covering $\tilde S\rightarrow S$ such that $\deg(\phi_{\tilde S})=27$ and $q(\tilde S)=1.$
\end{theorem}

To achieve this, we work with the equations of $Z,$ $S$ given in \cite{BK}, \cite{BY} to find equations for the curves that pullback to generators
of the canonical system of $\tilde Z,$ $\tilde S,$ and we show that their intersection is empty.
The calculations are very demanding and we had to find several workarounds in order to succeed.

Remark: Sai-Kee Yeung's proof \cite{Ye} for the case $d=36$ is not correct.
Recently, he has informed me that he has a new proof that is also based on Borisov-Keum equations.

The computations for the case $d=27$ are harder than the ones for $d=36.$
They require the computation of equations of some fibres of the Albanese fibration of the Cartwright-Steger surface $S.$
More precisely, we compute, over a finite field, the equations of the three fibres that are fixed by the $\mathbb Z/3$ action of $S.$
Then we show that they are smooth, which answers a question from Cartwright-Koziarz-Yeung \cite[Corollary 5.3, Remark 5.6]{CKY},
in particular it implies that the Albanese fibration of $S$ is stable.

All computations are implemented with the computer algebra system Magma \cite{BCP}, and can be found on arXiv:1903.03017 as ancillary files.

We use the symbol $\equiv$ for linear equivalence of divisors, the rest of the notation is standard in Algebraic Geometry.

\bigskip
\noindent{\bf Acknowledgements}

\noindent The author thanks Lev Borisov for a useful correspondence and for providing the equations of the two ball-quotient surfaces from \cite{BK}, \cite{BY}.

\noindent This research was supported by FCT (Portugal) under the project PTDC/MAT-GEO/2823/2014,
the fellowship SFRH/BPD/111131/2015 and by CMUP (UID/ MAT/00144/2019),
which is funded by FCT with national (MCTES) and European structural funds through the programs FEDER,
under the partnership agreement PT2020.

\section{Lift to rationals}

There is a classical method for computing a rational number $x$ from its values modulo a set of primes, by combining Chinese
remaindering with Farey sequences (see e.g. algorithm 2 in \cite{BDFP}). It works well provided the set of primes is big and none
of these is a 'bad prime'. We have implemented this algorithm with Magma, the usage is \verb|LiftToRationals(n,p)|, where
$p$ is a list of prime numbers and $n$ is a list containing the values of $x$ modulo $p.$ We use it in the computations below.

\section{The case $\deg(\phi)=36$}

In \cite{BK}, Borisov and Keum give the equations of a fake projective plane $Z$, embedded in $\m P^9$ by its bicanonical system.
It is known that this surface has an action of $\m Z/3$ such that the quotient $Y:=Z/(\m Z/3)$ is a surface with invariants $p_g=0$ and $K^2=3$,
and with singular set the union of $3$ ordinary cusps ($\mathrm A_2$ singularities). 

Let $G_Z$ and $G_Y$ be the groups such that $Z=\m B/G_Z$ and $Y=\m B/G_Y,$ where $\m B$ is the unit ball in $\m C^2.$
Computing the index $4$ subgroups of $G_Z,$ we see that there is a unique normal subgroup $G_{\tilde Z}$ of $G_Y$ such that
$G_Y / G_{\tilde Z}\cong (\m Z/2)^2\times \m Z/3.$ Let $\tilde Z:=\m B/G_{\tilde Z}.$
We have an abelian covering that factors as $$\tilde Z\xrightarrow{\text{$(\m Z/2)^2$}} Z\xrightarrow{\text{$\m Z/3$}} Y.$$
Since $G_{\tilde Z}$ is a subgroup of the fundamental group of $Z,$ the $(\m Z/2)^2$-covering is \'etale.
This gives $\chi(\tilde Z)=4$ and $K_{\tilde Z}^2=36.$
The maximal abelian quotient of $G_{\tilde Z}$ is a finite group, thus $q({\tilde Z})=0$ and then $p_g({\tilde Z})=3.$

Our goal is to show that the canonical map of $\tilde Z$ is of degree $36$ onto $\m P^2.$
This happens if and only if the canonical system of $\tilde Z$ is free from base points.
By \cite[Proposition 4.1]{Pa}, this system is generated by the pullback of three curves in $Y.$
Let $C_1,C_2,C_3$ be the corresponding curves in the fake projective plane $Z.$
Notice that $C_i$ is linearly equivalent to $K_Z$ up to $2$-torsion, thus $2C_i\equiv 2K_Z$ and then
$2C_i$ is a hyperplane section of $Z\subset\mathbb P^9.$
We will find the equations of these hyperplanes and verify that $C_1\cap C_2 \cap C_3=\emptyset,$
which implies that $|K_{\tilde Z}|$ is free from base points.

The curves $C_i$ are invariant for the $\m Z/3$ action.
Keeping the notation from \cite{BK}, let $\mathbb P^9=\mathbb P^9(U_0,\ldots,U_9)$ and define
$\m Z/3$-invariant sections
$$X_1:=U_1+U_2+U_3,\ \ X_2:=U_4+U_5+U_6,\ \ X_3:=U_7+U_8+U_9.$$
We need to search for hyperplane sections $H_i$ of $Z$ of the type
\begin{equation}\label{eqhyp}
a_0U_0+a_1X_1+a_2X_2+a_3X_3=0
\end{equation}
and such that $H_i=2C_i,$ $i=1,2,3.$
The strategy is to work over a finite field $\mathbb F_p$ and test all possible values of $a_0,\ldots,a_3$.
Then after finding a solution, repeat it for enough values of $p$, and finally use our Magma function \verb|LiftToRationals|
to obtain the solution over characteristic zero.
\\

\noindent{\bf Step 1.}\\
Let $C_1$ be the reduced subscheme of the scheme defined in \cite[Remark 2.2]{BK}, and let $H_1$ be the hyperplane section of $Z$ given by $U_0=0.$
We use the Magma function \verb|Difference| to show that $C_1=H_1-C_1,$ thus $H_1=2C_1.$
\\

\noindent{\bf Step 2.}\\
For each possibility for the coefficients $a_1, a_2, a_3$, we need to check if the hyperplane $H$ given by (\ref{eqhyp}) is not reduced. This is very time consuming, thus we test instead if $C_1\cap H$ is reduced or not. Notice that here we can remove $U_0$ from the equation of $H$, because $C_1$ is contained in the hyperplane $U_0=0$. Then we assume $a_1=1.$
Since the degree of $C_1$ is $18,$ we search only for the cases where the degree of the reduced subscheme of $C_1\cap H$ is at most $9.$
\\

\noindent{\bf Step 3.}\\
We compute this for several different values of the prime number $p,$ obtaining two solutions for each $p:$ $a_2=a_3=0$ and $a_2,a_3\ne 0$.
With such data we use our Magma function
\verb|LiftToRationals| and obtain the liftings
$$a_2=a_3=0\ \ {\rm and}\ \ a_2=\frac{1}{2}\left(\sqrt{-7}-3\right), a_3=\frac{1}{8}\left(\sqrt{-7}+5\right).$$

\noindent{\bf Step 4.}\\
For each of these two cases, we need now to test all hyperplanes $H$ of the type $a_0U_0+X_1+a_2X_2+a_3X_3=0$, running over all possible values of $a_0\in\mathbb F_p\left(\sqrt{-7}\right)$. In order to speed up computations, we take the hyperplane $H_{U_1}$ of $Z$ cut out by $U_1=0$ and test if $H_{U_1}\cap H$ is reduced or not.
Since the degree of $H_{U_1}$ is $36,$ we search for the cases where the degree of the reduced subscheme of $H_{U_1}\cap H$ is at most $18.$
\\

\noindent{\bf Step 5.}\\
We repeat for several different values of $p$ to obtain a list of pairs $p,a_0$. Then we use again the Magma function \verb|LiftToRationals|,
obtaining the hyperplanes
$$\left(1-\sqrt{-7}\right)U_0+4X_1=0$$ and
$$\left(-\sqrt{-7}-5\right)U_0+32X_1+\left(16\sqrt{-7}-48\right)X_2+\left(4\sqrt{-7}+20\right)X_3=0.$$
\\

\noindent{\bf Step 6.}\\
Let $H_2,$ $H_3$ be the corresponding hyperplane sections of the fake projective plane $Z,$ defined over $\m Q\left(\sqrt{-7}\right).$
We need to show the existence of curves $C_2, C_3\subset Z$ such that $H_2=2C_2, H_3=2C_3.$
Ideally we would just compute the reduced subscheme $R$ of $H_i$, $i=2,3,$ but our computer cannot finish this task.
Our workaround here is to find the system of quadrics through $R$, as follows.
We compute a zero-dimensional subscheme $v$ of $R$ with degree big enough such that every quadric that contains $v$ must also contain $R$.
Then we define the subscheme $C_i$ of $H_i$ cut out by these quadrics, and show that $H_i=2C_i$ (using the Magma function \verb|Difference|).
\\

\noindent{\bf Step 7.}\\
Consider the $2$-torsion divisors $t:=C_1-C_2$ and $t':=C_3-C_1.$ If $t\equiv t',$ then $2C_1\equiv C_2+C_3,$
which is impossible because there is no hyperplane through $C_2+C_3.$
This confirms that $t,t'$ generate the group $(\m Z/2)^2,$ that corresponds to the covering $\psi:\tilde Z\rightarrow Z.$

Finally we check that $C_1\cap C_2\cap C_3=\emptyset.$
Since the curves $\psi^*(C_i),$ $i=1,2,3,$ generate the canonical system of $\tilde Z,$
then the canonical map of $\tilde Z$ is of degree $36$ onto $\m P^2.$

\section{The case $q=1, \deg(\phi)=27$}

In \cite{BY}, Borisov and Yeung give the equations of the so-called Cartwright-Steger surface $S$, embedded in $\m P^9(U_0,\ldots,U_9)$
by its bicanonical system (we keep their notation).
It is known that this surface has an action of $\m Z/3$ such that the quotient $\bar X:=S/(\m Z/3)$ is a surface with
singular set the union of six ordinary cusps ($\mathrm A_2$ singularities) with three $\frac{1}{3}(1,1)$ singularities,
and whose smooth minimal model has invariants $p_g=1,$ $q=0$ and $K^2=2$.
Correspondingly there is a $\m Z/3$ Galois covering $$\psi:S\rightarrow\bar X.$$

Borisov and Yeung also give the equations of the unique effective canonical divisor of $S,$ it is the reduced subscheme of the hyperplane
of $S$ given by $U_0=0.$ We let $K_1$ be this curve.

It is known that the surface $\bar X$ contains a pencil of curves with three multiple fibres $F_i'=3D_i,$ $i=1,2,3,$
such that $F_1'$ contains the three $\frac{1}{3}(1,1)$ singularities, and $F_2',F_3'$ contain three cusps each.
One has $\psi^*(F_i')=3F_i,$ $i=1,2,3,$ where each $F_i$ is a fibre of the Albanese fibration of $S$.

Since two points in an elliptic curve move in a pencil, the same happens for $2F_1.$ We explicitly compute below the
pencil $|2F_1|,$ and show that it contains the divisor $F_2+F_3$
(this linear equivalence could be proved by using the fact that there is a unique elliptic curve with
an automorphism of order 3 that fixes points).
 This implies that $F_1-F_2\equiv F_3-F_1.$
Consider the $3$-torsion element $F_1-F_2$ and the corresponding \'etale $\m Z/3$ Galois covering $$\varphi:\tilde S\rightarrow S.$$

Let $G_{\tilde S}$ and $G_S$ be the groups such that ${\tilde S}=\m B/G_{\tilde S}$ and $S=\m B/G_S,$ where $\m B$ is the unit ball in $\m C^2.$
Computing all index $3$ subgroups of $G_S,$ we see that the maximal abelian quotient of $G_{\tilde S}$ is $\m Z/7\times\m Z^2$
or $\m Z^2$, thus $q({\tilde S})=1$ and then $p_g({\tilde S})=3.$

We want to show that the canonical map of $\tilde S$ is of degree $27$ onto $\m P^2.$
This happens if and only if the canonical system of $\tilde S$ is free from base points.
By \cite{Pa}, this system is generated by the pullbacks of three curves $K_1,K_2,K_3\subset S,$
with $K_1\equiv K_S.$
These are linearly equivalent up to $3$-torsion. We will compute $K_2, K_3$ as elements in
$$|K_1+F_1-F_2|,\ \ \ |K_1+F_1-F_3|,$$ respectively.
Finally we will verify that $K_1\cap K_2 \cap K_3=\emptyset,$
which implies that $|K_{\tilde S}|$ is free from base points.

The computation of $K_2$ and $K_3$ is very demanding, we have succeeded only working over finite fields $\m F_p$.
Fortunately, we got that $K_2+K_3$ is the hyperplane of $S$ with equation
$$U_7-2U_8-4U_9=0$$
for several different values of $p,$ which suggests that it remains unchanged over the rationals.

In the next section we show how to compute the equations of $K_2$ and $K_3,$ working over $\m F_p.$
Then we take the equations of $K_1$ and $K_2$ over the rationals, and do the necessary verifications.

\subsection{Computation of the hyperplane $K_2+K_3$}

Here we work over a finite field $\m F_p.$ We first compute the linear system $|K_1+F_1|$ (which is of dimension 19),
then the systems $|K_1+F_1-F_2|$ and $|K_1+F_1-F_3|$, giving the curves $K_2$ and $K_3$, respectively.
\\

\noindent{\bf Step 1.}\\
JongHae Keum \cite{Ke} shows that a fibre of the Albanese fibration of $S$ is numerically equivalent to $-E_1+5E_2,$
where $E_1,$ $E_2$ are certain irreducible curves. Lev Borisov has informed me that $E_1+E_2$ is the subscheme
of $S$ cut out by the hyperplane $\{U_1=0\}.$ Then Magma gives the prime components of this hyperplane, i.e. the equations of $E_1$ and $E_2.$

We use the Magma function \verb|IsLinearSystemNonEmpty| to compute the unique element in the linear system $|-E_1+5E_2|.$
This curve contains the three points of $S$ that correspond to the three $\frac{1}{3}(1,1)$ singularities of $S/(\m Z/3),$
thus it is the fibre $F_1.$
\\

\noindent{\bf Step 2.}\\
From the equations of $F_1,$ it is easy to give the defining equations of $2F_1,$ but we want an equation with the lowest possible degree and
not identically zero on $S$. We use the Magma function \verb|Divisor| to get a basis $B$ of the ideal of $2F_1$
(this takes several hours to finish).
We then choose one polynomial $g_9\in B$ of degree $9$, the lowest possible degree, and take the corresponding hypersurface $H$ of $S$.
\\

\noindent{\bf Step 3.}\\
Let $C$ be such that $H=2F_1+C$. We compute the basis of the ideal of the divisor $C,$
from where we take another degree $9$ polynomial $g_9'$ containing $C$
such that $g_9,g_9'$ generate the pencil $|2F_1|$ (after removing the base component $C$).
\\

\noindent{\bf Step 4.}\\
There is one element in this pencil containing six points that are fixed by the action of $\mathbb Z/3.$
After removing the base component $C,$ it must be the union of two Albanese
fibres, thus it is $F_2+F_3.$ In this way we obtain the equations of $F_2+F_3.$ 
\\

\noindent{\bf Step 5.}\\
Now we compute the divisor $K_1+F_1,$ and then we use the Magma function \verb|RiemannRochBasis| to compute a basis of its space of global sections.
This basis is generated by some rational functions, with numerators $N_i$, and with a common denominator. These are given on affine coordinates, so we take the projective closure.
\\

\noindent{\bf Step 6.}\\
Let $L$ be the linear system generated by the $N_i$.
We compute the unique element of $L$ that contains $K_1+F_1$ and take the corresponding curve in $S,$ say given by $N_1$.
Then we compute the intersection of $S$ with the base scheme of $L,$ which is $B=N_1-K_1-F_1$.
\\

\noindent{\bf Step 7.}\\
We want to compute the element of $L$ that contains the fibre $F_2$, but we don't have a factorization of the curve $F_2+F_3$, thus we use a workaround: the factorization of the zero-dimensional scheme $(F_2+F_3)\cap E_1$ contains two irreducible schemes of degree 39. We guess that one of these is in $F_2$ and the other is in $F_3$. Then we compute the element $N_2$ of $L$ through the first one (hence through $F_2$).
\\

\noindent{\bf Step 8.}\\
Finally we use the Magma function \verb|Complement| to compute the effective divisor $K_2:=N_2-B-F_2$ (which satisfies $K_2\equiv K_1+F_1-F_2$).
\\

\noindent{\bf Step 9.}\\
We repeat the above steps in order to get the curve $K_3$ ($\equiv K_1+F_1-F_3$).
\\

\noindent{\bf Step 10.}\\
Looking to the equations, we verify that $K_2+K_3$ is cut out on $S$ by the hyperplane
\begin{equation}\label{eqhypS}
U_7-2U_8-4U_9=0.
\end{equation}

\subsection{Linear equivalence of $3K_1,3K_2,3K_3$}

Here we work over the rational field.
\\

\noindent{\bf Step 11.}\\
We get the defining equations of the curves $K_2$ and $K_3$ by computing the prime components
of the hyperplane of $S$ given by (\ref{eqhypS}).
\\

\noindent{\bf Step 12.}\\
A straightforward computation gives the system $J_3$ of degree 3 hypersurfaces that contain the divisor $3K_1$.
Note that any element of $J_3$ is $\equiv 6K_1.$ 
\\

\noindent{\bf Step 13.}\\
For $i=2,3,$ we show that $3K_1+3K_i\in J_3,$ which implies $3K_i \equiv 3K_1.$
\\

\noindent{\bf Step 14.}\\
The fact $K_2+K_3\equiv 2K_1$ gives $K_2-K_1\equiv K_1-K_3,$ thus the curves $K_1, K_2, K_3$ pullback to linearly equivalent curves in $\tilde S.$

Finally we check that $K_1\cap K_2\cap K_3=\emptyset.$
Since the curves $\varphi^*(K_i),$ $i=1,2,3,$ generate the canonical system of $\tilde S,$
the canonical map of $\tilde S$ is of degree $27$ onto $\m P^2.$

\section{The $\mathbb Z/3$-invariant Albanese fibres}

Working over a finite field, the Magma function $\verb|IsLinearSystemNonEmpty|$ gives a unique element in each of the systems $|K_1+F_1-K_2|$ and $|K_1+F_1-K_3|,$ which are then the curves $F_2$ and $F_3.$ We want to compute the singular subset of the fibres $F_1,F_2,F_3.$ Since a direct computation is hard, we proceed as follows.

Consider the map $\rho:S\rightarrow\mathbb P^3$ given by $(U_0:\cdots:U_4).$
Let $X$ be the minimal resolution of the surface $\bar X:=S/(\mathbb Z/3).$ One can show that $\rho$ is the composition of the triple covering $S\rightarrow X$ with the bicanonical map of $X$. This bicanonical map is birational, hence $\left.\rho\right|_{F_i}$ is of degree 3. We check that the images $\rho(F_i)$ are smooth. This implies that the fibres $F_i$ can be singular at most at the $9$ points of $S$ that are fixed by the $\mathbb Z/3$ action. The computation says that this is not the case, thus $F_1,F_2,F_3$ are smooth.

\bibliography{References}

\

\

\noindent Carlos Rito
\vspace{0.1cm}
\\{\it Permanent address:}
\\ Universidade de Tr\'as-os-Montes e Alto Douro, UTAD
\\ Quinta de Prados
\\ 5000-801 Vila Real, Portugal
\\ www.utad.pt, {\tt crito@utad.pt}
\vspace{0.1cm}
\\{\it Temporary address:}
\\ Departamento de Matem\' atica
\\ Faculdade de Ci\^encias da Universidade do Porto
\\ Rua do Campo Alegre 687
\\ 4169-007 Porto, Portugal
\\ www.fc.up.pt, {\tt crito@fc.up.pt}\\

\end{document}